\input amstex
\documentstyle{amsppt}
  \catcode`@11



\catcode`@12

\let\olddate=\date
\def\date#1{\olddate#1\enddate}
\let\oldthanks\thanks
\def\thanks#1{\oldthanks#1\endthanks}
\let\oldabstract\abstract
\def\abstract#1{\oldabstract#1\endabstract}

\def\martinstricks{}   

\newtoks\fooline
\def\today{\number\day.\number\month.\number\year}
\fooline={{\tt\jobname}\hfill{\tenrm\relax\folio\hfill\today}}

\catcode`@11

\newtoks\oldoutput
\oldoutput=\the\output
\output={\setbox255\vbox to \ht255 
         {\unvbox255 %
\iffirstpage@\else 
\smallskip\line{\the\fooline}\vss\fi}\the\oldoutput}

\catcode`@12

\martinstricks

\def\demo#1{\smallskip\noindent{\smc #1}}
\def\enddemo{\medskip}

\define\Levy{\hbox{\rm L\'evy}}
\define\La{\Cal L}
\redefine\L{\Cal L}
\define\name#1{\vphantom{\underset\sim\to#1}
{\smash{\underset\sim\to#1}\vphantom{#1}}}
\def\Q{\name Q}
\topmatter
\title The hanf Numbers of Stationary Logic II:
\\Comparison with other Logics \endtitle
\author Saharon Shelah \endauthor
\affil The Hebrew University, Jerusalem Israel \\ Simon
Fraser University, Burnaby, B.C., Canada \\ Rutgers
University, New Brunswick, New Jersy, U.S.A. \endaffil
\date{ Aug 15, 1991 }
\abstract {We show the ordering of the Hanf number of
$\L_{\omega, \omega}(wo)$, (well ordering)
$\L^c_{\omega, \omega}$ (quantification on countable
sets), $\L_{\omega, \omega}(aa)$ (stationary logic) and
second order logic, has no more restraints provable in
$ZFC$ than previously known (those independence proofs
assume $CON(ZFC$ only).  We also get results on corresponding
logics for $\L_{\lambda, \mu}$.}
\thanks {The author would like to thank the BSF and NSREC
for partially supporting this work.  Publ number 211 }
\endtopmatter
\document
\hsize=4.5in
\magnification=1200
\baselineskip=18pt
\vfill\pagebreak

\subhead \S 0 Introduction \endsubhead

The stationary logic, denoted by $\La(aa)$ was
introduced by Shelah [Sh 43].  Barwise, Kaufman and
Makkai [BKM] make a comprehensive research on it, proving
for it the parallel of the good properties of $\L(Q)$. 
There has been much interest in this logic, being both
manageable and strong, see [K] and [Sh 199].

Later some properties indicating its afinity to second
order logic were discovered.  It is easy to see that
coutable cofinality logic is a sublogic of $\L(aa)$. 
By [Sh 199],  for pairs  $\varphi,\psi$ of formulas in
$\L_{\omega,\omega}(Q^{cf}_{\aleph_0})$, satisfying
$\vdash\varphi\to\psi$ there is an interpolant in
$\L(a,a)$.  By Kaufman and Shelah [KfSh 150],
for models of power $>\aleph_1$, we can express in
$\L_{\omega,\omega}(aa)$ quantification on countable
sets.  Our main conclusion is (on the logics see Def 1.1
or the abstract, on $h$, The Hanf Numbers, see 1.2)
\bigskip
\proclaim {0.1  Theorem} The only restiction on the Hanf
numbers of $\L_{\omega,\omega}(wo),
\L^c_{\omega,\omega},
\L_{\omega,\omega}(aa),\L^{II}_{\omega,\omega}$

 are:

$h\bigl(\L_{\omega,\omega}(wo)\bigr) \leq
h\bigl(\L^c_{\omega,\omega}\bigr) \leq
h\bigl(\L_{\omega,\omega}(aa)\bigr) \leq
h\bigl(\L^{II}_{\omega,\omega},h)$

 $h(\La^c_{\omega,\omega}) <h(\La^{II}_{\omega,\omega}$.  
\endproclaim
\bigskip
\demo{Proof} See {\bf 2.1} (neccessity), {\bf 2.2, 2.4,
2.5} and {\bf 3.3} (all six possibilities are consistent).

The independence results are proved assuming $CON(ZFC)$
only and the results are generalized to
$\L_{{\lambda^+},\omega}$.  We do not always remember to
write down the inequalities of the form $\L_{\lambda,
\omega}(Q_1) <\L_{\mu,\omega}(Q_2)$.  For some of the
results when we generalize them to
$\L_{{\lambda^+},\omega}$ or
$\L_{\lambda, \kappa}$ we need a stronger	
hypothesis.  The proofs of the results on 
$h(\L_1) \leq h(\La_2) $ give
really stronger information: we can interpret
$\L_1 $ in $\L_2$, usually here
by using extra predicates, i.e.,every formula in
$\La_1$ is equivalent to a formula in
$\triangle (\La_2)$; remember $\triangle (\La_2)$ is defined by: $\theta \in
\triangle (\La_2) (\tau)$ is represented
by $(\theta_1, \theta_2),\, \theta_e \in
\La_2(\tau_e), \tau_1\cap\tau_2=\tau,\,
M\vDash \theta$ iff $ M$ can be expanded to a model of
$\theta_1$ iff $M$ cannot be expanded to a model of
$\theta_2$ ( so the requirement on $(\theta_1,
\theta_2)$ is strong).  Note that this has two
interpretation:  one in which we allow $\tau_1, \tau_2$
to have new sorts hence new elements, the other in which we do not allow it. We
use an intermediate course, we allow this but the  number of new elements are
the power set of the old.  But for   $\La^c_{\omega,\omega} \leq \La_{\omega,
\omega}(aa)$, for models of power $\lambda=\lambda^{\aleph_0}$ we do not need
new elements.

We thank Matt Kaufman for discussions on this subject.
\enddemo
\bigskip
\definition{Notation}
 Let cardinals be denoted by
$\lambda, \kappa,\mu,\chi$

Ordinals are denoted by $\alpha, \beta, \gamma, \xi,
\zeta,i, \jmath$. $\delta$ is a limit ordinal.

Let $H(\lambda)$ be the family of sets whose transitive
closure has cardinality $<\lambda$ (so for $\lambda$
regular it is a model of $ZFC~$, i.e., $ZFC^-$ except the
power set axiom: and for a strong limit a model of $ZC$.

Let L\'evy$(\lambda, \kappa)=\{f:\, f$ a function
from some $\alpha < \lambda$ into $\kappa\}$

\qquad L\'evy $(\lambda, {<}\kappa)=\{f:\,f$ a
partial function from $\lambda \times \kappa$ to $\kappa,
|Dom f| <\lambda,  f(\alpha,\beta) <1+\beta\}$.
\enddefinition

\definition{Notation on Logics}: $\La$ will be
a logic, $\tau$ a vocabulary (i.e., set of predicates and
fuction symbols, always with a fixed arity, usually finite).  We assume
that $\La(\tau)$ is a set of formulas, each
with $ < O c_1(\La)$ free variables and
$<Oc(\La)$  predicates and function
symbols; $\La(\tau)$ is closed under first order
operations, substitutions and relativizations and
$\La(\tau) $ is a set (with $\tau$ and the
 the family of variables sets)

Two formulas are isomorphic if some mapping from the set of
predicates, function symbols and free variables of one onto
those of another is one-to-one and map one formula to the
other.

We are assuming that up to isomorphism there is a set of
$\La$-formulas, this number is denoted by
$|\La|$.

Let $\La_1\subseteq \La_2$ mean
$\La_1 (\tau) \subseteq \La_2(\tau)$
for every vocabulary $\tau$.
\enddefinition
\vfill\pagebreak

\subhead \S1 \, Preliminaries \endsubhead
\proclaim{1.1 Definition}
\roster
\item $\La_{\lambda,\kappa}$ is the logic in which
$\wedge_{i\in I}(|I|<\lambda)$ and
$(\exists x_0,\dots,x_i\dots)_{i \in J}
(|J|<\kappa)$ are allowed, with
$Oc_1(\La_{\lambda,\kappa})=\kappa $ (so
$\La_{\omega,\omega}$ is first order logic)
\item  For a logic $\La, \La(wo)$ extends $\La$ by
allowing the quantifier $(wo\, x, y)\varphi(x,y)$ saying
$\langle\{x:\exists\varphi(x,y)\}, \, \varphi (x,y)\rangle$ is well
ordering
\item For a logic $\La, \La^c=\La(\exists^c)$ extends
$\La$ by allowing a monadic predicate as free variable
and the quantifier $(\exists^c X)\varphi(X)$ saying ``there
is a countable set $X$ such that $\varphi(X)"$
\item For a logic $\La, \La(aa)$ extends $\La$ by
allowing  monadic predicates as free variables and the
quantifiers $(aaX)\varphi(X)$ saying that the collection of
countable $X$ satisfying $\varphi$ contains a closed
unbounded family of countable subsets of the model
\item For a Logic $\La, \La^{II}=\La(\exists^{II}) $
extends $\La$ by allowing binary predicates as free
variables and the quantifiers $\exists R\varphi(R)$saying
there is a two-place relation $R$ on the model
satisfying $R$
\item For $Q \in \{\exists^c, aa, \exists^{II}\},
\La'(Q)$ is defined similarly allowing a string

$(Q x_1\dots Q x_i \dots)_{i<\alpha},
|\alpha|< 0c_1(\La)$
\item Let $\La^c=\La(\exists^c),\La^{wo}=\La(wo),
\La^{II}=\La(\exists^{II}), \La^{aa}=\La(aa)$
\endroster
\endproclaim
\bigskip
\definition{1.2 Definition}
\roster
\item    For a sentence $\psi$ Let
$h(\psi)=sup\{|M|^+:M\models \psi\}$

(so it is a cardinal (or infinity)) and it is the first
$\lambda$ such that $\psi$ has no model $\geq \lambda)$
\item   For a theory $T, h(T)=h(\wedge_{\psi \in T} \psi)$
\item   For a logic $\La$ let $
h(\La)=\text{sup}\{h(\psi):h(\psi) < \infty,\psi \in
\La(\tau)$ for some vocabulary $\tau\}$
\item    For a logic $\La$  and cardinal $\lambda$ let
$h(\La,\lambda)=\text{sup}\{h(\psi): $ for
some vocabulary $\tau$ of power $<\lambda, \psi \in
\La(\tau), h(\psi)<\infty\}$
\item   For a logic $\La$ and cardinal $\lambda$
let $hth(\La, \lambda)=\text{sup}\{h(T): $ for some
vocabulary $\tau$ of power $<\lambda, T \subseteq \La
(\tau), h(\tau) <\infty\}$

$hth(\La)=H(\La, \infty)$
\endroster
\enddefinition
\bigskip
\proclaim{1.3 Claim}
\roster 
\item  for every $\psi \in \La$ for some $\varphi \in
\La, \, [h(\psi)<\infty \rightarrow h(\psi) < h(\varphi) <
\infty]$ \item $h(\La)$ is strong limit
\item If $\La$ is closed under $\wedge_{\alpha
<\alpha_0}$ for $\alpha_0<\lambda$ then $cf[h(\La)] \geq
\lambda$
\item If the number of sentences in $\La$ (up to
isomorphism) is $\leq \lambda$ then $cf[h(\La)] \leq
\lambda$
\endroster
\endproclaim
\bigskip
\proclaim{1.4 Lemma}  assume $\La$ is a logic $\subseteq
\La^{II}_{\omega,\omega}$ and there is a function $f$
from Card to Card such that:
\roster
\item"(a)" $f$ is definable in
$\La^{II}_{\omega,\omega}$, i.e., the class of two
sorted models $\langle \kappa, f(\kappa) \rangle$ is definable by some
sentence of $\La^{II}_{\omega,\omega}$ or even just
\item"(a)$^-$" For some
$\lambda^*<h(\La^{II}_{\omega,\omega})$ and $\varphi^* \in
\La^{II}_{\omega,\omega}$ for $\kappa, \mu, \geq
\lambda^*,\kappa < h(\La^{II}_{\omega,\omega})$
we have $\langle \kappa,\mu \rangle \vDash\varphi^*$ iff $\mu=f(\kappa)$
\item"(b)" If $\psi \in \La$ has a model of power $\geq
\kappa$ {\it then} $\psi$ has a model $M,
\kappa \leq ||M|| \leq f(\kappa)$
\item"(c)" $\La$ is definable in
$\La^{II}_{\omega,\omega}$ i.e., the class $\{(\psi,
\tau,M):\,\psi\in \La(\tau), \, Ma \, \tau$-model,
$M\vDash \psi\}$ is definable by a sentence in $\La^{II}_{\omega,\omega}$
\item"(d)" For $\mu < h(\La), f(\mu) <h(\La) $

\noindent\underbar{Then} $h(\La)  <h(\La^{II}_{\omega,\omega})$
\endroster
\endproclaim
\bigskip
\demo{Proof} Easy.  Let $\psi_0 \in
\La^{II}_{\omega,\omega}$ be such that
$\lambda^*<h(\psi_0) < \infty$, where $ \lambda^*, \varphi^*$
are as in (a)$^-$.  We can assume $h(\psi_0)< h(\La)$
(otherwise the conclusion is trivial).  Let $\psi \in
\La^{II}_{\omega,\omega}$ say that for some $\lambda,
\mu_0$:
\roster
\item"(i)" the model $M$ is isomorphic to some
$(H(\lambda),\in),\lambda$ strong limit,
\item"(ii)" for every $\kappa <\lambda,
M\models(\exists\mu \geq \kappa)[\psi_0$ has a model of
cardinality $\mu] \vee (\exists \mu \geq \kappa)
[\langle\kappa,\mu\rangle\models \varphi^*]$
\item"(iii)" $\mu_0<\lambda , \psi_0$ has a model of
power whose cardinality is in the interval $\in
(\mu_0,\lambda)$
\item"(iv)" for every $\kappa < \lambda, \kappa\geq
\mu_0$, there is $\theta \in \La$ which has a model of
cardinality in the interval $(\kappa,\lambda)$, but for
some $\kappa' \in (\kappa,\lambda)$ has no model of
cardinality in the interval $(\kappa',\lambda)$
\endroster
Now $(H(h(\La)),\in)$ is a model of $\psi$ and it has no
models of larger cardinality.\vfill $\boxdot_{1.4}$

We can prove similarly:

\proclaim{1.5 Lemma}  Suppose $\La_1,\La_2$ are logics
and there is $f:\text{Card}\to$ Card such that
\roster
\item"(a)"  for some $\lambda^* <h(\La_2)$ and $\varphi^*\in
\La$ for $\kappa,\mu \geq\lambda^*$ we have: $\langle\kappa,\mu\rangle \models
\varphi^*$ iff $\mu=f(\kappa)$
\item"(b)"  if $\psi \in \La_1$ has a model of cardinality $ \geq \kappa$ {\it
then} $\psi$ has a model $M, \kappa \leq ||M||\leq f(\kappa)$ 
\item"(c)"  $\L_1$ is definable in
$\La_2$ just in the following weaker sense: for $K_1=\{(\psi, \tau):\psi \in
\La_1(\tau)\},\,K_2=\{M,\psi,\tau):M\models  \psi,\,\psi \in \La_1(\tau)\}$
there are $\psi_e \in \La_2$.

$(\forall x)[x \in K_e \Leftrightarrow$ for some
$\lambda$, some expansion of $(H(\lambda),\in,x)$
satisfies $\psi_e]$ and for every $x \,\{\lambda:
\text{some expansion of}\, (H(\lambda),\in, x)$ satisfies
$\psi_e\}$ is a bounded family of cardinals
\item"(d)"  For $\mu <h(\La_1), f(\mu) < h(\La_2)$
\item"(e)" $||\La|| <h(\L_2), \La_{\omega,\omega}^{II}
\subseteq \La_2$  

\noindent\underbar{Then} $h(\La_1) <h(\La_2)$
\endroster
\endproclaim
\bigskip
\remark{Remark}  Of course if  1.5 is
hypothesis holds for $\La_1$ (and $\La_2$) then the
conclusion holds for $\La'_1,\La'_2$ whenever $\La'_1 \subseteq \La_1$ and
$\La_2 \subseteq \La'_2$.
\endremark

\proclaim{1.7 Lemma}\roster
\item If $ M\vDash \psi, \psi \in \La^{wo}$ then this is
preserved by any forcing, this holds even for $\psi \in
\La^{wo}_{\infty,\omega}$
\item If $M \vDash \psi, \psi \in
\La_{\omega,\omega}^{aa}$ then this is preserved by any
$\aleph_1$-complete forcing this holds even for $\psi \in
\La^{aa}_{\infty,\omega}$
 \item  If $M \vDash \psi, \psi\in
\La_{\omega,\omega}^c$ this is preseved by forcing not adding new countabale
subsets of $|M|$ (this holds even for $\psi \in
\La_{\infty ,\omega_1}$)
\item  If $M\vDash \psi, \psi \in \La_{\infty,\lambda},
\lambda$ regular, then this is preserved by forcing by
$P$ where P does not add sequences of ordinals of lenght
$<\lambda$.  If $P$ is  $\aleph_1$-complete this holds
for $\psi \in \L_{\infty,\lambda}^{aa}$.
\item Suppose $V_1,V_2$ are models of set theory (with
the same ordinals), $V_1 \subseteq V_2$, and letting
$\lambda=h(\La)^{V_1}$ where $\La$ is
$\La^{wo}_{\omega,\omega}$ or $\L^c_{\mu,\omega}$ or
$\La^c_{\mu,\omega}$, (just a suitable downward Lowenheim
Skolem theorem is needed).

{\it If} $\{A \subseteq \lambda : A$ bounded,
$A \in V_1\}=\{A\subseteq \lambda: A$ bounded, $A\in
V_2\}$ {\it then} $h((\La)^{V_1}=h(\La)^{V_2}$.
\endroster
\endproclaim
\demo{Proof}  Left to the reader.
\enddemo
 
\subhead \S 2 Independence for $\La^c_{\omega,\omega},
\La_{\omega\omega}^{II}$ \endsubhead

In this section we shall deal with the indepedence of the
cases where
$h(\La_{\omega,\omega}^{wo})=h(\La^c_{\omega,\omega})$.

\proclaim{2.1 Lemma}
\roster
\item   For any logic $\La: h(\La(wo))\leq h(\La^c)
\leq h(\La(aa))\leq h(\La^{II})$
\item   For any logic $\La$ we have $h(\La^c_{\omega\omega})
<h(\La^{II}_{\omega,\omega})$ 
\item   For any logic $\La$ we have $h(\La^c_{{\lambda^+},\omega}) <
h(\La^c_{\lambda^+,\omega})$,  moreover:

if $\lambda < h(\La^{II}_{\mu,\omega})$ {\it then}
$h(\La^c_{\lambda^+,\omega}) < h(\La^{II}_{\mu,\omega})$
\endroster
\endproclaim
\bigskip
\demo{Proof}
\roster
\item   By Kaufman and Shelah [KfSh 150, Theorem 4.1]; only
$\La=\La_{\omega,\omega}$ is discussed there, but it
makes no difference, the non trivial part is $h(\La^c)
\leq h(\La^{aa})$;

\item  See [KfSh 150];
\item  Use {\bf 1.5} for the function
$f:f(\kappa)=(\kappa^{\aleph_0})^+$
\endroster
\enddemo
\bigskip
\proclaim{2.2 \ Lemma}
\roster
\item  If $V=L$ then
$h(\La^{wo}_{\omega,\omega})=h(\La^c_{\omega,\omega})
<h(\La^{aa}_{\omega,\omega})=h(\La^{II}_{\omega,\omega})$
\item  If $V=L$, then for any logic $\La,\quad
h(\La^{wo})=h(\La^c) \leq h(\La^{aa})=h(\La^{II})$.
\endroster
\endproclaim
\bigskip
\demo{Proof}
\roster
\item See [KfSh 150]
\item Same proof.
\endroster
\enddemo
\bigskip 
\proclaim{2.3 \ Fact} For a regular cardinal $\lambda$
and $\psi\in \La_{\lambda,\lambda}^{aa}$ the following
are equivalent:
\roster
\item"(i)"  for every $\mu$ large enough
$\Vdash_{\Levy(\lambda,\mu)} ``\psi \text{has a model of
power}\,\lambda$''
\item"(ii)"  for some $\lambda$-complete forcing $Q$ we have:
$\Vdash_Q ``\psi \,\text{has a model of power}\, \geq
\lambda$''.
\endroster
\endproclaim
\bigskip
\demo{Proof}  Easy; (i)$\Rightarrow $ (ii): as
L\'evy$(\lambda,\mu)$ is a $\lambda$-complete forcing notion, (i) is a
particular case of (ii). (ii)$\Rightarrow$(i) let $Q$ be a
$\lambda$-complete forcing notion such that  $\Vdash_Q ``\psi$ has a model of
cardinality $\geq \lambda"$.  Let $\mu$ be such that $\mu >|Q|, \Vdash_Q
``\psi$ has a model of cardinality $\geq \lambda$ but
$\leq \mu"$ and $\mu=\mu^\lambda$.  In $(V^Q)^{\text{L\'evy}(\lambda,\mu)}
 \psi$ has a model of cardinality
$\lambda$ by {\bf 1.7(4)}. 

 But  $(V^Q)^{\text{L\'evy}(\lambda,\mu)}$ is $V^{\text{L\'evy}(\lambda,\mu)}$. 
(see e.g. [Kun]).
\enddemo
\bigskip
\definition{2.3A Notation}  Let $\mu_0[\psi,\lambda]$ be the
first  cardinal $\mu$ satisfying {\bf 2.3(i)}, if one exists, and
$\lambda$ otherwise.
\enddefinition
\bigskip
\proclaim{ 2.4 \ Lemma}
\roster
\item  In some forcing extension of $L,
h(\La^{wo}_{\omega,\omega})=h(\La^c_{\omega,\omega}) <
h(\La^{aa}_{\omega,\omega}) < h(\L^{II})$
\item Moreover for $\lambda <h(\La^{II}_{\mu, \kappa})$,
we have $h(\L^{aa}_{\lambda,\lambda}) <h(\La
^{II}_{\mu,\kappa})$
\endroster
\endproclaim
\bigskip
\remark{2.4A Remark}  If we want to have: $\lambda
< h(\La ^{aa}_{\mu,\lambda})\Rightarrow
h(\La^c_{\lambda,\omega}) < h(\La^{aa}_{\mu,\omega})$, we
should define
$\lambda_{i+1}=h(\La^c_{{\mu^+_i},\omega})^+$.
\endremark
\bigskip
\demo{Proof}  Start with $V=L$.  Let $\psi^*\in
\La^{aa}_{\omega,\omega}$ a sentence such that
$h(\La^c_{\omega,\omega}) <h(\psi^*) < \infty$ be chosen
later.  Let $\lambda_0 > h(\psi^*)$ be regular,
$\lambda_\circ < h(\La^{II}_{\omega,\omega})$.  We define
an iterated forcing $ \langle P_i,
\Q_j:i \leq\infty, \jmath <
\infty\rangle$ and cardinals $\lambda_i$ such that:
\roster
\item"(a)"  the iteration is with set support (so
$P_\infty$ is a class forcing)
\item"(b)"  $\lambda_i$ is
regular cardinal
\item"(c)"  $\lambda_i \geq \sum_{\jmath
<i} \lambda_\jmath$, and $\lambda_i$ is the
first regular cardinal $\geq\sum_{j<i}(\lambda_j+\mu_j)^+$ (when $i>0)$
\item"(d)"  $Q_i (\in\,V^{P_i})$ is
$\lambda_i$-complete
\item"(e)"  Let $\{\psi^i_\alpha:\alpha <
\lambda_i\}$ be the set of all
$\La^{aa}_{{\lambda_i},{\lambda_i}}$ sentences
(up to isomorphism) in $V^{P_i}$.
\endroster
We define in $V^{P_i}, Q_i$ to be $\text{L\'evy}
(\lambda_i, \mu_i)$ where $ \mu_i$ is
the successor of $sup\{\mu_0[\psi,\lambda_i]^{V^{P_i}}:\psi
\in \La_{{\lambda_i},{\lambda_i}}\}$ and so
$\lambda_{i+1}=\mu_i^+$.

Our model is $V^{P_\infty}$.  Clearly the $\lambda_i$
are not collapsed (as well as limits of $\lambda_i$
and $\chi<\lambda_0$) and other successor cardinals $\geq
\lambda_0$ are collapsed.  So in $V^{P_\infty}$, for
regular $\chi \geq \lambda_0$, if $\psi \in
\La^{aa}_{\chi,\chi}$ has a model of cardinality $\geq\chi$ then it has a model
of cardinality $\chi$.  As clearly
$h(\La^{II}_{\omega,\omega}) > \lambda_0$, we get by {\bf
1.4} $h(\La^{aa}) <h(\La^{II})$ (as well as ({\bf 2})).

By the Lowenheim Skolem theorem, using {\bf 1.7(5)} for
$\psi \in \La^{wo}_{\omega, \omega}$ or $\psi \in \La
^c_{\omega,\omega}, h(\psi)$ does not change (being
$\infty$ or $<\lambda_0$) hence (in $V^{P\infty}$)\ 
$h(\La^{wo}_{\omega,\omega})=h(\La^{wo}_{\omega,\omega})^V;
h(\La^c_{\omega,\omega})=h(\La^c_{\omega,\omega})^V$.
Hence (in
$V^{P_\infty})\,h(\La_{\omega,\omega}^{wo})=h(\La^c_{\omega,
\omega})$ as this holds in $L$.

We still have to choose $\psi^*\in
\La^{aa}_{\omega,\omega}$ and prove that in $V^{P_\infty}$ we have
$h(\La^c_{\omega,\omega})<h(\La^{a,a}_{\omega,\omega})$.  There is $\psi^*\in
\La^{aa}_{\omega,\omega}, L\vDash "h(\L^c_{\omega,\omega}) < h(\psi^*) <
\infty$" (by {\bf 2.2}).

Clearly for any such $\psi^*, V^{P_\infty} \vDash
``h(\La^c_{\omega,\omega}) < h(\psi^*)"$ (as no new subset
of $h(\psi^*)$ is added), but we need also
$V^{P_\infty}\vDash`` h(\psi^*) <\infty$"; but checking the
sentneces produced in [KfSh 150] proof of {\bf Theorem 4.3}
(for proving $L \vDash h(\La^{aa})=h(\La^{II})$),
they are like that.  So $V^{P_\infty}\vDash
"h(\La^c_{\omega,\omega}) <h(\La^{aa}_{\omega,\omega})"$.
\enddemo
\bigskip
\proclaim{2.5\ Lemma}
\roster 
\item In some forcing expresion of $L$ we have
$h(\La^{wo}_{\omega,\omega})=h(\La^c_{\omega,\omega})=
h(\La^{aa}_{\omega,\omega})< h(\La^{II}_{\omega,\omega})$
\item In fact for any logic $\La$ we have
$h(\La^{wo})=h(\La^c)=h(\La^{aa})$
\item  For $\lambda <h(\L^{II}_{\mu,\kappa})$ {\it then},
$h(\La^{aa}_{\lambda,\lambda})=h(\La^{II}_{\mu,\kappa})$.
\endroster
\endproclaim
\bigskip
\demo{Proof}  We start with $V=L$. We define a (full set
support) iteration, $\bar Q=\langle P_i,\Q_i:i$
an ordinal $\rangle (\Q_i-a\, P_i$ name) and
cardinals $\lambda_i$ such that 
\roster
\item"(a)" $ \lambda_i$ is regular $\geq \aleph_1 +|P_i|$ for  $i$ limit
$\lambda_i=(\sum_{j<i}\lambda_i)^+$
 \item"(b)" $\Q_i$ is  $\lambda_i$-complete
\item"(c)"  if $i $ is even, $G_i \subseteq
P_i$ the generic set (remember
$\Q_i\in V^{P_i})$ then let the set of
elements of $P_i$ be listed as
$\{p_\alpha^i:\alpha < \lambda_i\}$, and
$\Q_i$ will be the product of the L\'evy collapses of
$\aleph_{{\lambda_i \omega}+4\alpha+2+m}$ to
$\aleph_{{\lambda_i\omega}+4\alpha+1+m}$ 
for  $\alpha <\lambda_i$ such
that: \,$[p_\alpha^i \in G_i \Rightarrow m =0]$
and $[p^i_\alpha \not\in G_i \Rightarrow m=1]$.
Let
$\lambda_{i+1}=\aleph_{{\lambda_i\omega}+{\lambda_i}+1}$
\item"(d)"  if $i$ is odd, let
$\{\psi_\alpha^i: \alpha <\lambda_i\}$ list
all sentneces of
$\L^{aa}_{{\lambda_i},{\lambda_i}}$ in a
 rich enough vocabulary of cardinality $\lambda_i$).  For each $\alpha$ if there
is a $\lambda_i$-complete forcing notion $Q$ (which
is a set) and (in $V^{P_i}) \Vdash_Q$ ``there is a
model of $\psi_\alpha^i$ of cardinality $\geq
\lambda_i$'' then let $\mu^i_\alpha$ be such
that  $\Vdash_{\text{L\'evy}({\lambda_i},{\mu^i_\alpha}})
``\psi^i_\alpha \text{ has a model of cardinality}
\lambda_i"$; otherwise $\mu^i_\alpha=\lambda_i$.

Note that $\mu^i_\alpha$ exists by {\bf 2.3}.

Let $Q_i=\text{L\'evy}(\lambda_i,
<\lambda_{i +1})$ where
$\lambda_{i+1}=(\lambda_i
+\sum{\alpha<\lambda_i} \mu^i_\alpha)^{++}$.
\endroster
Let $G_\infty \subseteq P_\infty$ be generic over $V$ and
$V[G_\infty]$ be our model.  Note in $V[G_\infty]$, 

{\bf (*)} $[i\, \text{odd} \Rightarrow
\lambda_{i+1}=\lambda_i^+]$

\qquad $[i\, \text{even} \Rightarrow
\lambda_{i+1}=\lambda_i^{+({\lambda_i}\omega+1)}]$

\qquad $[i\, \text{limit} \rightarrow \lambda_i
=(\sum_{\jmath<i} \lambda_\jmath)^+]$.

For  $\lambda=\lambda_{2\jmath+1}$, if $\psi\in
\La^{aa}_{\lambda,\lambda}$ has a model of cardinality
$\geq\lambda$ then it has a model of cardinality
$\lambda$ (by {\bf 2.3 + 1.7(4)}).  By {\bf (*)} we deduce
that $V^{P_\infty} \vDash$ ``if $\psi \in
\L_{\lambda,\chi}^{aa}$ has a model of cardinality
$>\lambda$ then it has a model $M, \lambda < ||M|| < \aleph
_{\lambda^+}"$.

So {\bf 1.5} is applicable to show
$h(\La^{aa}_{\omega,\omega}) <h(\La^{II}_{\omega,\omega})$
(and by {\bf 1.6} and {\bf 1.7}) also {\bf 2.5(3)} holds.

Why $h(\La^{wo}_{\omega,\omega})=\La^{a,a}_{\omega,\omega}?$.  Let
$\psi^*$ describe  $(L_\lambda \in, G_\infty \cap \cup_{i<\delta} P_i)$.

If $M \vDash \psi^*$, then for some $\alpha$ and
$G,M\cong(L_\alpha, \in,G)$, so without loss of
generality equality holds.  Now if $\lambda <|\alpha|, M
\vDash ``\lambda$ is a [regular] cardinal of
$L"$ {\it iff} $\lambda$ is a [regular] cardinal
of $\La$.  Also we know that for every ordinal $\zeta$, if
in $L,\lambda_{2i} \leq \aleph_\zeta <
\lambda_{2i+1},  \zeta$ divisible by four
 then forcing by
$P_\infty$ collapses at most one of the cardinals
$\aleph_{\zeta+1},\aleph_{\zeta+2},\aleph_{\zeta
+3},\aleph_{\zeta+4}$ of $L$; if $\lambda_{2i} \omega \leq\zeta <\lambda_{2i}
\omega+\lambda_{2i}$ then exactly one.

We assume $\psi^*$ say so, and so when
$\aleph_{\zeta+4}^L$  $\leq |\alpha|$ the answer in
$M$ to the  question ``which of
$\aleph_{\zeta+1},\aleph_{\zeta+2}, \aleph_{\zeta+3},
\aleph_{\zeta+4}$ is collapsed'' is the right one.  So
when $\lambda_{2i+1} <|\alpha|$, we can in $M$
reconstruct $G_\infty \cap P_{2i}$ (see choice of
$Q_{2i}$).

But $V^{P_\infty} \vDash ``\lambda_{2i+1} \leq
\aleph_{\lambda_{2i}(\omega+1)+1}$ and
$\lambda_{2i+2}=(\lambda_{2i+1})^+$ and for limit
$\delta$ we have
$\lambda_\delta=(\sum_{i<\delta}\lambda_i)^+$"

The rest is as in [KfSh 150] proof of {\bf 4.3}
\vfill\pagebreak

\subhead \S 3 {\bf $h(\La^{wo}_{\lambda,\omega})$} is O.K. but for
{\bf $h(\La_{{\aleph_3},\omega})$} large cardinals are needed
and sufficient \endsubhead

  In section {\bf 2} we deal with the three cases for
which
$h(\La^{wo}_{\omega,\omega})=h(\La^c_{\omega,\omega})$. 
Here we deal with the three cases where
$h(\L^{wo}_{\omega,\omega}) <h(\L^c_{\omega,\omega})$. 
The new part is {\bf Lemma 3.2}, and then, in {\bf 3.3}
we get the desired conclusion.  For dealing with
$\La_{{\lambda^+},\omega}$ we do not assume $CON(ZFC)$
alone, we assume the existence of a class of large
cardinals (weaker than measurability).  By {\bf 3.4} at
least if $\lambda\geq \aleph_3+ (2^{\aleph_0})^+$,
something of this sort is necessary.
\bigskip
\proclaim{3.1\ Fact}:  The following are equivalent for
$\psi \in \L^{wo}_{\omega,\omega}$ or even $\psi \in
\L^{wo}_{\infty,\omega}$:
\roster
\item"(i)" for every $\mu$ large enough $\Vdash_{\hbox{L\'evy}
{(\aleph_0},<\mu)} ``h((\psi) =\infty$"
\item"(ii)"  for some (set) forcing notion $P$ we have $\Vdash_P
``h(\psi)=\infty$''.
\endroster
\endproclaim
\bigskip
\demo{Proof}  similar to the proof of {\bf 2.3}
\bigskip
\remark{3.1A\ Notation}  Let the first $\mu$ satisfying (i)
be $\mu_1(\psi)$ (and $\aleph_0$ if there is no such
$\mu$).
\endremark
\bigskip
\proclaim{3.2\ Lemma}  $(V=L)$. 
\roster \item For some (set) forcing
notion $P$
$$\Vdash_P\,``h(\La^{wo}_{\omega,\omega})
<h(\La^c_{\omega,\omega})$$
and this is preserved by
$h(\La^{wo}_{\omega,\omega})^+$-complete forcing".
\item In (1) we can use L\'evy$(\aleph_0<\mu)$ for some $\mu>cf\mu=\aleph_0$
\item We can use instead Cohen$(\mu)=\{f:f $ a finite function from $\mu$ to
$(0,1)\}$. So cardinals are not collapsed
\endroster 
\endproclaim
\bigskip
\demo{Proof} 1) Let \ $\mu^*=sup\{\mu_1(\psi):\psi\in
\L^{wo}_{\omega,\omega}\}$

We now define a finite support iteration $\langle P_i,
\Q_n :n<\omega\rangle$  and 
$\mu_n$ as follows:

\qquad $\mu_0=\mu^*$

\qquad $Q_0=\text{L\'evy}(\aleph_0,\mu_0)$ 

\qquad for $n\geq 0,\mu_{n+1}$ is
$h(\La^{wo}_{\omega,\omega})^{V^{P_n}}$

\qquad $Q_n=\text{L\'evy}(\aleph_0,\mu_n)$.

Let $\mu=(\sum \mu_n)$.  Note that $P_\omega$ satisfies
the $\mu^+-c.c$.

Now $V^{P_\omega}$ is our model.  Note

{\bf (*)}  $V^{P_\omega} \vDash G.C.H.\,+\aleph_1=\mu^+,
\text{and}\, V=L[\Bbb R, <]$ for any well $<$
ordering of $\Bbb R$.

Note that in $\goth B=(\omega \cup \Cal P
(\omega))^{V^{P_\omega}};o,+,\times,\in)$ we
can define by first order formulas (representing
ordinals by well ordering of $\omega$):
\roster
\item"(a)"  $\cup_n \mu_n$ (maximal countable
ordinal which is a cardinal in $L_{\mu^+}$
\item"(b)"  $L_{\mu^+} $ hence $\langle
\mu_n:n<\omega\rangle$ (by induction remembering the Lowenheim Sholem theorem)
hence the iteration (really we can omit this as
$P_\omega$ is just L\'evy ($\aleph_0,<\mu)$)
\item"(c)"  the set $\Bbb R^-=^{def}\{r \in \Bbb R:$
for some $n$, and $G \subseteq P_n$
generic over $V, r \in V[G]\}$.

And for $r \in \Bbb R^-$
\item"(d)" $H_r=\{\psi \in \La^{wo}_{\omega,\omega}:
L[r] \vDash h(\psi) < \infty\}$ as it is equal
to

$\{\psi \in \La^{wo}_{\omega,\omega}: L[r]
\vDash h(\psi) <\cup_n\mu_n\}$.
\endroster
[Note that $P_n 's$ are homogeneous, hence
$h(\psi)$ does not depend on $G\subseteq P_n$] 

So by {\bf 3.1} and the choice of $\mu_0$, we
can define in that model $\goth B$

$H^*=\{\psi \in
\La^{wo}_{\omega,\omega}:h(\psi)^{V^{P_\omega}}<
\infty\}$

[How ? it is $\cap\{H_r:r\in \Bbb
R^-\}$,remembering {\bf 3.1}]

Let
$\lambda=h(\La^{wo}_{\omega,\omega})$ (in $
V^{P_\omega})$.

Now we define a sentence $\varphi \in
\La^c_{\omega,\omega}$: it just describes
$(H(\lambda),\in)$: it says
\roster
\item"(i)" enough axioms of $ZFC$ holds
\item"(ii)" every countable bounded set of ordinals is
represented
\item"(iii)"  on every infinite cardinal $\alpha$
there is a model $M_\alpha$ with universe $\alpha$
satisfying some  $\psi \in H^*$ (which we have shown is
definable in any model $M$  of $\varphi)$
\endroster

 So we have proved the first assertion from {\bf 3.2}.  Now $\lambda$-complete
forcing, preserve trivially $``h(\psi) \geq\mu"$ as it preserves
satisfaction for $\La^{wo}_{\omega,\omega}$.  It preserves
$``h(\psi)<\infty"$ as this is equivalent to
$``h(\psi)<\lambda"$, the forcing adds no new model power
$<\lambda$, and Lowenheim Skolem Theorem finishes the
argument.

2)  We have proved it in the proof of (1)

3)  A similar proof, replacing $\mu_1 (\psi)$ by $\mu'_1=$ first $\mu$ such
that $\Vdash_{Cohen(\mu)}``h|\psi|=\infty"$
if there is one $\aleph_0$
otherwise.\vfill $\square_{3.2}$
 \proclaim{3.3 Conclusion}  for some forcing extensions of $\L$:
\roster
\item $h(\La^{wo}_{\omega,\omega})
<h(\La^c_{\omega,\omega})
<h(\La^{aa}_{\omega,\omega})=h(\La^{II}_{\omega,\omega})$
\item $h(\La^{wo}_{\omega,\omega})
<h(\La^c_{\omega,\omega})=h(\La^{aa}_{\omega,\omega})
<h(\La^{II}_{\omega,\omega})$
\item $h(\La^{wo}_{\omega,\omega}) <h(\La^c_{\omega,\omega})
<h(\La^{aa}_{\omega,\omega}) <h(\La^{II}_{\omega,\omega})$.
\endroster
\endproclaim
\bigskip
\demo{Proof}: Combine {\bf 3.2} with \S 2.
\enddemo
\bigskip
\proclaim{3.4 Claim} ($\neg 0^\#)$: For $\lambda \geq
\aleph_3 + (2^{\aleph_0})^+$ we have
$h(\La^{wo}_{\lambda,\omega})=h(\La^c_{\lambda,\omega})$.
\endproclaim
\bigskip
\remark{Remark}:  The logics are essentially equivalent.
\endremark
\bigskip
\demo{Proof}  If $\psi \in \La^{wo}_{\lambda,\omega}$
says $M$ is, for some $\alpha, (L_\alpha[A],\in)$ (up
to isomorphism), $\alpha >2^{\aleph_0}, A \subseteq
2^{\aleph_0}$, every subset of $\omega$ is in
$L_{(2^{\aleph_0})} [A]$, and $\alpha \geq \omega_2$,
and $\{\delta < \aleph_2:cf\delta=\aleph_0$ in
$L_{\omega_2}[A]\}=\{\delta
<\aleph_2:cf\delta=\aleph_0\}$ then by Jensen's covering
lemma $[\beta<|\alpha|\Rightarrow$ every countable subset of
$\beta$ is represented in the model \hfill$\square_{3.4}$
\enddemo
\bigskip
\proclaim{3.5 Claim}  Suppose that:
\roster
\item"(*)"  for every $\chi$ for some $\mu,
\mu\rightarrow(\omega_1)^{<\omega}_\chi$ or even just
\item"(**)"  for every $\chi$ for some
$\mu,\mu\rightarrow_{BG}(c)^{<\omega}_\chi$, which
means:  for every $f:[\mu]^{<\omega}\rightarrow\chi$ for
some $\langle\gamma_n: n<\omega\rangle$ for every
$\alpha<\omega_1$, for some $Y\subseteq \mu$, $Y$ has
order type $\alpha$ and $\wedge_n 
(\forall w\in [Y]^n)[\gamma_n=f(w)]$ .

 {\it Then} for every
$\lambda, h(\La^{wo}_{\lambda^+,\omega})$. 
\endroster
\endproclaim
\bigskip
\remark{3.5A Remark} 
\roster \item The property (**) was discoverd by
Baumgartner and Galvin [BG] such
that:

$\mu_{\rightarrow_{BG}}(c)^{<\omega}_\chi$ {\it iff}
$\mu \geq h(\La^{wo}_{\chi^+,\omega})$.
\item See [KfSh 150, 4.2] (for $\lambda=\omega$)
\endroster
\endremark
\bigskip
\demo{Proof}  There is a sentence $\psi \in
\La^c_{\lambda,\omega}$ such that for $\chi \leq
\mu$: there is a model $M, ||M||=\mu, |P^M|=\lambda$, {\it
iff} $(\forall\alpha < \mu ) \alpha\rightarrow
_{BG}(c)^{<\omega}_\chi$.\hfill $\square_{3.5}$
\bigskip
On $K=K^V$ (the core model of V) see Dodd and Jensen [DJ].
\enddemo
\proclaim{3.6 Claim}  Suppose $V=K$, and (**)  ( from
{\bf 3.5} ),then
\roster
\item for every $\lambda$ we have
$h(\L^{wo}_{\lambda^+,\omega})  <
h(\L^c_{\lambda^+,\omega}) <
h(\L^{aa}_{\lambda^+,\omega})=
h(\La^{II}_{\lambda^+,\omega})$
\item  for every $ \La, h(\La^{aa})=h(\La^{II})$.
\endroster
\endproclaim
\bigskip
\demo{Proof} 1)  First inequality by the observation
above, the second inequlity follows from last equality
{\bf Th 2.1}, last equality see {\bf( 2)} 

(note: if $cf\delta>\aleph_0$ in $\La^{aa}_{\omega,\omega}$ we can say for
$A\subseteq\delta$ whether $\{\alpha<\delta:cf\alpha=\aleph_0, \alpha \in A\}$ is
a stationary subset of $\delta$).

2)  As in [KfSh 150]
\enddemo
\bigskip
\remark{3.7 Observation}  There is $\psi \in
\La^c_{\omega,\omega}$ such that $M\vDash \psi$ {\it iff}
$M$ is isomorphic to $K_\alpha$ for some $\alpha$.

It is known see (see [BG], [DJ])
\endremark
\bigskip
\definition{3.8 Fact}  If in $V$ there are, e.g.,
measurable cardinals in Card, then $K\vDash$ (**).
\enddefinition
\bigskip
\proclaim{3.9 Claim}:  Suppose $V=K$ and (**) holds.

For some forcing extension $V[G_\infty]$ of $V,
V[G_\infty]\vDash (**)$ and for every $\lambda,
h(\La^{wo}_{\lambda^+,\omega}) <h(\La^c_{\lambda^+,\omega})
<h(\La^{aa}_{\lambda^+,\omega})
<h(\La^{II}_{\lambda^+,\omega})$.
\endproclaim
\bigskip
\demo{Proof}  Similar ot {\bf 2.4(1)} except that we want
to preserve (**).  We define by induction on $\alpha$ an
iterated forcing, $\langle P_i,Q_j \leq
\alpha,\jmath <\alpha\rangle$ with set support and cardinals
$\lambda_i$ increasing such that:
\roster
\item"(i)" $\lambda_0=\aleph_2$
\item"(ii)" $\lambda_\delta=(\sum_{i<\delta}\lambda_i
+|P_\delta|)^+$
\item"(iii)"  if $\lambda_i, P_i$ are defined,
let $\mu_i$ be $\lambda_i^+
+\cup\{\mu_0[\psi,\lambda_i]:\psi \in
\La_{\lambda_i,\lambda_i}^{aa}\}$.

$Q_i=$ L\'evy $(\lambda_i^+,\mu_i^+)
(\text{in}\, V^{P_i})$ and $\lambda_{i+1}$ is
minimal such that

$\lambda_{i+1}\rightarrow_{BG}(c)^{<\omega}_{\mu^{++}_i}$ and $\lambda_{i+1} \leq
h(\La^c_{{\mu^+_i}},\omega)$.
\endroster
We leave the rest to the reader.
\enddemo
\bigskip
\proclaim{3.10 Claim}  Suppose $V=K$ and (**) holds.

For some forcing extension $V[G_\infty]$ of $V,
V[G_\infty]\vDash(**)$ (hence the conclusion of {\bf 3.7})
and for every $\lambda$
$$h(\La^{wo}_{\lambda^+,\omega})
<h(\La^c_{\lambda^+,\omega})=
h(\La^{aa}_{\lambda^+,\omega})
<h(\La^{II}_{\lambda^+,\omega})$$
\endproclaim
\bigskip
\demo{Proof}  Combine the proofs of {\bf 3.9} and {\bf
2.5}.
\enddemo
\bigskip
\subhead \S 4 Lowering consistency strength \endsubhead
\bigskip
 We present here some alternative
proofs with lower consistency strength than in \S 3. Specifically {\bf 4.1,
4.3} and {\bf 3.2(3)} justfy the restriction $\lambda \geq
\aleph_3+(2^{\aleph_0})^+$ in {\bf 3.4].

\proclaim {4.1 Lemma} Let $V=L$.  Then there is a forcing notion $P\in L$, not
adding reals, such that for $G\subseteq P$ generic over
$V$, in $V[G]$:
\roster
\item"a)"  $h(\La^{wo}_{\omega_1,\omega})<h(\La^c_{\omega_1,\omega})$
\item"b)"  No $\aleph_1$-complete forcing notion (or
even forcing notion satisfying the $\Bbb I$-condition
$\Bbb I$ a set of $\aleph^{V[G]}_2$-complete ideals
from $L$ changes the truth value of "$h(\psi) <\infty$"
for $\psi \in L^{wo}_{\omega_1,\omega}$
 \item"c)"  There is a sentence $\psi \in
\La^c_{\omega,\omega}$ whose class of models of power
$\geq \aleph_2$ is just $\{L_\alpha[G]:\alpha
\geq\aleph_2\}$ ( and note $P\in
L_{\aleph_2}{[V[G]]})$
\item"d)"  $h(\La^c_{{\lambda^+},\omega})
<h(\La^{aa}_{\lambda^+,\omega})=h(\La^{II}_
{\lambda^+,\omega})$
\endroster
\endproclaim
\bigskip
\remark{4.1A Remark}  In the proof below, coding generic
sets by the decision which $L$-cardinals are collapsed is
replaced here by ``which $L$-regular cardinal have in $V$
cardinality $\aleph_0$ and which cardinality $\aleph_1$
\endremark
\bigskip
\demo{Proof}  Let $\Bbb I(\mu,\kappa)$ be, e.g. the calss
of filters $D$ which are $\lambda$-complete over some
$\lambda$ (this in $V$), where
$\mu\leq\lambda<\kappa, |\cup D| <\kappa$

We define by induction on $n,
\alpha_n,\beta_n,\lambda_{i,j},\mu_{i,j}, \langle P_i,
\Q_\jmath:i\leq\alpha_n,\jmath
<\alpha_n\rangle$ and $f_n$ such that
\roster
\item"(A)" $\alpha_0=0,\alpha_{n+1} >\alpha_n$
\item"(B)"
$\langle P_i,\Q_j,\mu_\jmath:i\leq
\alpha_n,\jmath <\alpha_n\rangle$ is an  $RCS$
iteration suitable for $ x_{\alpha_n}=\langle\Bbb I_{i,j}, \lambda_{i,j},
\mu_{i,j}^i, i<j\leq \alpha_n, i$ not strongly inaccessible $\rangle$.

See [Sh-b Ch.XI] or [Sh-f (Ch XI)] particularly Def. 6.1
\item"(C)" $f_n$ is a one-to-one function from $P_i$
onto some ordinal $\beta_n$, extending $\cup_{e<n} f_e$.

$G_\alpha$ will denote a generic subset of $P_\alpha$.

For $n=0$ there is nothing to do.

For $n+1$, note that forcing by $P_{\alpha_n}$ does not
add new reals.  So
$(\La^{wo}_{{\omega_1},\omega})^V=
(\La^{wo}_{{\omega_1},\omega})^{V[G_{\alpha_n}]}$ and let
$\{\psi_i:i <\omega_1\}$ be a lsit of the
sentences (up to isomorphism).
\endroster
We now (i.e for defining $\alpha_{n+1}$ etc.) define by
induction on $\zeta <\omega_1,
Q_{{\alpha_n}+\zeta}, x_{{\alpha_n}+\zeta+1}$ as follows:
\roster
\item"(a)" $\langle P_i, \Q_\jmath:i \leq
\alpha_n +\zeta\rangle$ is $x_{{\alpha_n}+\zeta+1}$-suitable $RCS$ iteration
\item"(b)" If there is $\Q_{{\alpha_n}+\zeta}$,
a $P_{{\alpha_n}+\zeta}$-name of a forcing notion
sattisfying the $\Bbb I
((|P_{{\alpha_n}+\zeta}|+sup\{\lambda_{i,\jmath}:i<j\leq\alpha_n+\zeta\})^+,
\kappa$
-condition
for some $\kappa$ {\it then}  $\Vdash_{{P_{\alpha_n+\zeta+1}\ast \Q}}
\psi_\zeta$ has arbitrarily large models then $\Q_{{\alpha_n}+\zeta}$
is like that, otherwise it is, e.g., L\'evy $(\aleph_1,2^{\aleph_1})$.
\endroster
Next let
$\mu_\zeta=h(\La^{wo}_{{\omega_1},\omega})^
{V[G_{{\alpha_n}+{\omega_1}}]},
Q_{{\alpha_n}+{\omega_1}}=$ L\'evy
$(\aleph_1,\mu^+_\zeta)$.  Now (where $<,>$ is
Godel's pairing function on ordinals), let in
$V[G_{{\alpha_n}+{\omega_1}} ]: A_n=|\{\langle f_n(p),
f_n(q)\rangle:p,q \in P_{\alpha_n} \vDash p\leq q$ and $p\not=q\}
\cup \{\langle f_n(p), f_n(p) \rangle:q\in
G_{{\alpha_n}+\omega}\}$ and let
$\gamma_n=sup\{\langle f_n(p),f_n(q)\rangle:p,q\in P_{{\alpha_n}+{\omega_1}}\}$. 
Now we define $Q_{{\alpha_n}+{\omega_1}+i}$ by induction on
$i\leq \gamma_n$:
\bigskip
 $Q_{{\alpha_n}+{\omega_1}}$ \quad is
L\'evy $(\aleph_1,\aleph_2)^{V[G_{{\alpha_n}+{\omega_1}}]}$,

 $Q_{{\alpha_n}+{\omega_1}+1+2i+1}$
\quad is L\'evy $(\aleph_1,
\aleph_2[V[G_{{\alpha_n}+{\omega_1}+1+2i+1}]]$,

$Q_{{\alpha_n}+{\omega_1}+1+2i}$
\quad is Namba forcing (of
$V[G_{{\alpha_n}+{\omega_1}+1+i}])$ if
$i \in A_n$ and  L\'evy
$(\aleph_1,\aleph_\zeta)^{V(n,i)}$ where
$V(n,i)=V[G_{{\alpha_n}+{\omega_1}+1+2i}]$
if $i \not\in A_n$.

Now let
$\alpha_{n+1}=\alpha_n+\omega_1+2\gamma_n,
\lambda_{n+1}=|P_{{\alpha_n}+{\omega_1}+2{\gamma_n}}|$,
and define $f_{n+1}$.

We leave the rest to the reader \hfill $\square_{4.1}$
\enddemo
\bigskip
\proclaim{4.2 Conclusion} 
\roster \item We can do the forcing from {\bf
2.4, 2.5} to the universe we got in {\bf 4.1}
getting corresponding results (for
$\La_{{\omega_1},\omega}(Q)$'s, with $CH$ and
$G.C.H)$: so we need $CON(ZFC)$ only.
\item  the same holds for {\bf 4.3} for the
$\La_{{\omega_2},\omega}(Q)$'s (so we use
$CON(ZFC+$ ``the class of ordinals in Mahlo'')
only).
\endroster
\endproclaim
\bigskip
\proclaim{4.3 Lemma}  Suppose $V=L$, (for
simplicity) and $\infty$ is a Mahlo cardinal (i.e.,
every closed unbounded class of cardinals has a
regular member).  {\it then} there is an
inaccessible cardinal $\lambda$ and a forcing notion
$ P\subseteq H(\lambda)$, such that :
\roster 
\item"(a)" $ P$ satisfies the $\lambda$-c.c., does
not add reals and collapse every $\mu \in
(\aleph_1,\lambda): and \Vdash_P \,"G.C.H.\lambda$ is
$\aleph_2 "$ and $|P|=\lambda$
\item"(b)"  $h(\La^{wo}_{{\omega_2},\omega})<h(\La^c_{\omega,\omega})$
\item"(c)" there is a sentence $\psi \in
\La^c_{\omega,\omega}$ whose class of models of power
$\geq \aleph_2$ is just suitable expansions of
$\{L_\alpha[G]:\alpha\geq \aleph_2\}$.
\endroster
\endproclaim
\bigskip
\demo{Proof}  Like {\bf 4.1}, but instead of induction on
$n<\omega$ we do induction on $\gamma <\infty$, and in
the induction only we first do the coding $(Q_{{\alpha_n}+
{\omega_1}+i},i <\gamma)$ (so that for $c$), we
say that for some club of $C$ of $\omega_2$, for $\delta
\in C$ we are coding the set of sentence in
$\L_{{|\delta|}^+}[G \cap P_\delta]$.
\enddemo
\bigskip

 Do we really need the large cardinal hypothesis in
{\bf 4.3} (and so in {\bf 4.2(2)})?
\proclaim{4.4 Claim}
Suppose $0^\#\not\in V$ and $\aleph_2^V$ is a successor
cardinal in $L$ and $2^{\aleph_0}=\aleph_1$ {\it then}
for some sentence $\psi \in \La^{wo}_{{\omega_2},\omega}$,
its models are exactly suitable expansions of
$(L_\alpha, \Cal P_{<{\aleph_1}}(\alpha))$, where $\alpha$
is the last $L$-cardinal $<\aleph^V_2$.

Hence
$h(\La^{wo}_{{\omega_2},\omega})=
h(\La^c_{{\omega_2},\omega})$.
\endproclaim
\bigskip
\demo{Proof}  Should be clear
\enddemo
\bigskip
\remark{4.4 Concluding Remarks}: Still we do not settle the exact consistency
strength.  In fact e.g. if $\aleph_2^V$ is the first $L$-inaccessible, we can
still prove the last sentence of {\bf 4.4}.

For $h(\La^{wo}_{{\omega_2},\omega})<h(\La^c_{{\omega_1},\omega})$ with
$2^{\aleph_0}=\aleph_1$ we can generalize {\bf Lemma 4.3} to this case (using
[Sh-f, XV]). 

  Also there is a gap in consistency strength in {\bf \S 3} for
$\lambda > \aleph_3 +(2^{\aleph_0})^+$.

It is not hard to show that if
$\lambda\geq\aleph_2+2^{\aleph_0},cf\lambda>\aleph_0$  and for some
$A\subseteq \lambda$ does not exists, then
$h(\La^{wo}_{{\lambda^+},\omega})=h(\La^c_{{\lambda^+},\omega})$
\endremark
\vfill\pagebreak
 

\Refs 
\ref \key [BKM] J.Barwise, M.Kaufman and M.Makkai,
\paper "Satationary Logic"
\jour Annals of Math Logic \vol13 \year 1978 \pages 171--224
\endref
\ref \key [BKM]
\paper A correction to "Stationary Logic"
\jour Annals of Math Logic
\vol 20 \year 1981 \pages 231--232
\endref
\ref \key [BG] J. Baumgartner and F. Galvin
\jour Annals of Math Logic
\endref
\ref \key  [DJ]  Dodd and Jenssen 
\paper The core model
\endref
\ref \key [K] M.Kaufman
\paper Model Theoretic Logics, J. Barwise and S. Feferman \vol 166 (6.12)
\endref
\ref \key [Kun] K. Kunen
\paper Set theory, An Introduction to Independence Proofs
\jour Studies in Logic and the Foundation of Math
\vol 102
\endref
\ref \key [KfSh 150] M. Kaufman and S. Shelah
\paper The Hanf Number of Stationary Logic, Part I
\jour Notre Dame J. of Formal Logic \vol 27 \year 1986 \pages 111--123
\endref
\ref \key [Sh 43] S. Shelah
\paper Generalized Quantifiers and Compact Logics
\jour Trans. Amer. Math Sci.
\vol 204
\year 1975
\pages 342--364
\endref
\ref \by [Sh b] S. Shelah
\paper Proper Forcing
\jour Springer Verlag, Lectures notes
\year 1982
\endref
\ref \by [Sh 199] S. Shelah
\paper Remarks in Abstract Model Theory
\jour Annals of Pure and Applied Logic
\issue 1985
\endref
\ref 
\by [Sh f]
\paper Proper and improper forcing \jour Springer Verlog, in prepint
\endref
\ref
\by [V] J. Vaananen
\paper On the Hanf numbers of unbounded logic
\endref
\endRefs
\enddocument